\documentclass[a4paper,12pt]{article}
\usepackage{amssymb}
\usepackage{color,enumitem,graphicx}
\usepackage{subfigure}
\usepackage{float}

\usepackage[colorlinks=true,urlcolor=blue,
citecolor=red,linkcolor=blue,linktocpage,pdfpagelabels,
bookmarksnumbered,bookmarksopen]{hyperref}

\newcommand{\R}{\mathbb{R}}

\newcommand{\cuad}{{\sqcap\kern-.68em\sqcup}}

\newtheorem{theorem}{Theorem}[section]
\newtheorem{proposition}{Proposition}[section]

\newtheorem{lemma}{Lemma}[section]

\newtheorem{remark}{Remark}[section]
\newcommand{\bremark}{\begin{remark} \em}
\newcommand{\eremark}{\end{remark} }

 \usepackage{cite}
\begin{document}

\begin{center}{\bf  \large    Positive supersolutions  for the Lane-Emden
system with inverse-square potentials  }\medskip
%%%%%%%%%%%%%%%%%%%%%%%%%%%%%%%%%%%%%%%%%%%%%%%%%%%%%%%%%%%%%%%%%%%%%%
%%%%%%%%%%%%%%%%%%%%%%%%%%%%%%%%%%%%%%%%%%%%%%%%%%%%%%%%%%%%%%%%%%%%%%
\bigskip\medskip

 {\small
 Huyuan Chen\,\footnote{{\tt chenhuyuan@yeah.net}}
}
 \medskip

{\small   Department of Mathematics, Jiangxi Normal University,\\
Nanchang, Jiangxi 330022, PR China
 } \\[3mm]

  {\small
 Vicen\c tiu D. R\u adulescu\,\footnote{{\tt radulescu@inf.ucv.ro}}
}
 \medskip

{\small  Faculty of Applied Mathematics, AGH University of Science and Technology,\\ 30-059 Krak\'ow, Poland
 } \\[3mm]

{\small  Department of Mathematics, University of Craiova,\\
          200585 Craiova, Romania
 } \\[3mm]

 {\small
Binlin Zhang\,\footnote{{\tt zhangbinlin2012@163.com}}
}
 \medskip

{\small   College of Mathematics and Systems Science, \\  Shandong University of Science and Technology,\\
Qingdao, Shandong 266590, PR China
 } \\[6mm]

 \begin{abstract}
In this paper,  we study the nonexistence of positive supersolutions for the following Lane-Emden system with  inverse-square potentials
 \begin{equation}\label{0}
 \left\{
\begin{array}{lll}
 -\Delta u+\frac{\mu_1}{|x|^2} u= v^p \quad {\rm in}\ \, \Omega\setminus\{0\},\\[2mm]
  -\Delta v+\frac{\mu_2}{|x|^2} v=  u^q \quad {\rm in}\ \, \Omega\setminus\{0\}
 \end{array}
 \right.
 \end{equation}
for suitable  $p,q>0$,    $\mu_1,\mu_2\geq -(N-2)^2/4$, where $\Omega$ is a smooth bounded domain containing the origin in $\R^N$ with $N\geq 3$.  Precisely, we provide
sharp supercritical regions of $(p,q)$ for the nonexistence of positive supersolutions to system (\ref{0})
in the cases  $-(N-2)^2/4\leq \mu_1,\mu_2<0$ and $-(N-2)^2/4\leq \mu_1<0\leq \mu_2$.  Due to the negative coefficients  $\mu_1,\mu_2$  of the inverse-square potentials,  an initial blowing-up at the origin could be derived and an iteration procedure could be applied in the supercritical case to improve the blowing-up rate until the nonlinearities are not admissible in some weighted $L^1$ spaces.  In the subcritical case, we prove the existence of positive supersolutions for system (\ref{s 1.1})
by specific radially symmetric functions.

\end{abstract}

\end{center}
%\tableofcontents \vspace{1mm}
  \vspace{1mm}


\vspace{2mm}

\setcounter{equation}{0}
\section{Introduction and main results}

Let the integer $N\ge 3$, $\mu\ge\mu_0:=-\frac{(N-2)^2}{4}$ and the Hardy operator be defined $\mathcal{L}_\mu=-\Delta+\mu|x|^{-2}$. Our concern of this paper is to consider  nonexistence of positive supersolutions for Lane-Emden system
\begin{equation}\label{s 1.1}
  \arraycolsep=1pt\left\{
\begin{array}{lll}
\mathcal{L}_{\mu_1}  u= v^p \quad &{\rm in}\quad \Omega\setminus\{0\},\\[2mm]
  \mathcal{L}_{\mu_2}  v= u^q \quad &{\rm in}\quad \Omega\setminus\{0\},
%\ u,\,v\geq 0\quad &{\rm on}\quad \partial\Omega,
 \end{array}
 \right.
 \end{equation}
where   $\mu_1,\,\mu_2\geq \mu_0$, $p,\,q>0$,  $\Omega$ is a bounded $C^2$ domain containing the origin  in $\mathbb{ R}^N$.

It is known that Hardy operator is related to the classical Hardy inequality  stated as following: for a smooth bounded domain $\mathcal{O}$ in $\R^N$ containing the origin, there holds
$$
\int_{\mathcal{O}} |\nabla u|^2 dx \ge c_{N}  \int_{\mathcal{O}}   |x|^{-2} |u|^2 dx ,\quad \forall\,u\in H_0^1(\mathcal{O})
$$
with the best constant $c_N=\frac{(N-2)^2}{4}$. The qualitative properties of Hardy inequality and its improved versions have been studied extensively, see for example \cite{ACR,BM,Fi,GP}, motivated by great applications in  the study of semilinear elliptic Hardy problem by variational method, see \cite{DD,D,FF}  and the reference therein.  Due to the inverse-square potentials,  the related semilinear elliptic and parabolic equations appear various peculiar phenomena and attract great attentions, such as  parabolic equations \cite{CM1,PV,VZ},  singular solutions of Hardy problem with absorption nonlinearity \cite{CV,CC,C,GV}, Hardy problems with source nonlinearities  \cite{BDT,CZ,F, FPZ}.  %More references on the isolated singularities see \cite{BV,V0,V1} .

When $\mu\ge\mu_0$,   the homogeneous problem
\begin{equation}\label{eq 1.1}
  \mathcal{L}_\mu u= 0\quad{\rm in}\quad  \R^N\setminus \{0\}
 \end{equation}
has   radially symmetric  solutions with the formula that
\begin{equation}\label{1.2}
\Phi_\mu(x)=\arraycolsep=1pt\left\{
\begin{array}{lll}
|x|^{\tau_-(\mu)}\quad&{\rm if}\quad \mu>\mu_0\\
|x|^{\tau_-(\mu)}(-\ln|x|) \quad&{\rm if}\quad \mu=\mu_0
\end{array}\right.
\quad{\rm and}\quad \Gamma_\mu(x)=|x|^{ \tau_+(\mu)},
\end{equation}
where
\begin{equation}\label{1.1}
 \tau_\pm(\mu)=- (N-2)/2\,\pm\sqrt{\mu-\mu_0}
\end{equation}
 are two roots of $\mu-\tau(\tau+N-2)=0$.  Note that $\tau_+(\mu)<0$ if $\mu\in [\mu_0,\, 0)$,  and any positive solution of
\begin{equation}\label{eq 1.1}
  \arraycolsep=1pt\left\{
\begin{array}{lll}
\mathcal{L}_{\mu}  u= u^p \quad &{\rm in}\quad \Omega\setminus\{0\},\\[1.5mm]
 \quad \ \,  u=0 \quad &{\rm on}\quad \partial\Omega
%\ u,\,v\geq 0\quad &{\rm on}\quad \partial\Omega,
 \end{array}
 \right.
 \end{equation}
 blows up at least like $\Gamma_\mu$ at the origin.  Inversely,  for the existence of positive solutions, the nonlinearity $u^p$ must have an admissible singularity in some weighted $L^1$ space.  Inspired by this observation,  \cite{BDT,D} show that letting $p^*_\mu=1+\frac{2}{-\tau_+(\mu)}$ for $\mu\in [\mu_0,\, 0)$, semilinear problem (\ref{eq 1.1})
 has no any  positive solution for $p\geq p^*_\mu$. This nonexistence of  (\ref{eq 1.1}) is
 specific to Hardy problem with $\mu\in[\mu_0,0)$.\smallskip

 Our interest in this paper is  to obtain the nonexistence of positive  supersolutions of (\ref{s 1.1}) when at least one parameter of $\mu_1,\mu_2$ is
 in $[\mu_0,0)$. Here a functions' pair  $(u,v)$ is a positive  supersolution of Lane-Emden system (\ref{s 1.1}), if $(u,v)$  satisfies  the inequalities
    \begin{equation}\label{s 1.1-def}
\mathcal{L}_{\mu_1}  u(x)\geq  v(x)^p\quad{\rm and}\quad
  \mathcal{L}_{\mu_2}  v(x)\geq u(x)^q \quad  {\rm for\ any \ } x\in  \Omega\setminus\{0\}.\end{equation}

 The nonexistence results state as follows.
 \begin{theorem}\label{teo 1}
 Assume that $\mu_0\leq \mu_1<0\leq \mu_2$ and $p,q>0$ satisfy one of the following assumptions:
$$
 \arraycolsep=1pt
\begin{array}{lll}
   (i) \quad  q\geq \frac{N+\tau_+(\mu_2)}{-\tau_+(\mu_1)}; \qquad\qquad\qquad\qquad\qquad\qquad\qquad\qquad\qquad\qquad\qquad\qquad\\[6mm]
   (ii)  \quad    \frac{2}{-\tau_+(\mu_1)}  <q<   \frac{N+\tau_+(\mu_2)}{-\tau_+(\mu_1)}\quad  {\rm and}\quad
  \left\{\begin{array}{lll} \tau_+(\mu_1) (pq-1)+2p+2<  0 \ \ {\rm if }\ \, \mu_1>\mu_0,\\[2mm]
 \tau_+(\mu_1) (pq-1)+2p+2\leq 0\ \ {\rm if}\ \, \mu_1=\mu_0,
 \end{array}\right.
 \qquad
  \end{array}
$$
 then system (\ref{s 1.1}) admits no positive supersolutions.

\end{theorem}

\begin{figure}[H]
    \centering
   \subfigure{
  \begin{minipage}{65mm}
   \includegraphics[scale=0.48]{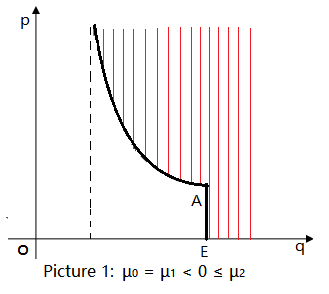}

  \end{minipage}
  }
    \subfigure{
    \begin{minipage}{60mm}
  \includegraphics[scale=0.50]{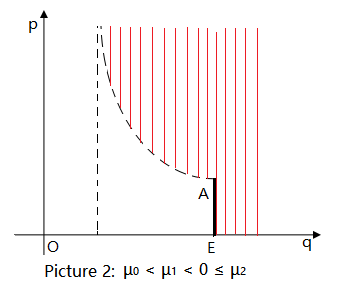}

  \end{minipage}
  }

  {\footnotesize   The red region of $(p,q)$ in Picture 1 and 2 represents the nonexistence of positive solution of system (\ref{s 1.1}).  }
 \end{figure}

\begin{theorem}\label{teo 2}
 Assume that $\mu_0\leq \mu_1,\,\mu_2<0$ and $p,q>0$ satisfy one of the following assumptions:
 $$
 \arraycolsep=1pt
\begin{array}{lll}
   (i) \qquad \ \ \ \ {\rm either}\quad\ p\geq \frac{N+\tau_+(\mu_1)}{-\tau_+(\mu_2)}\quad{\rm or}\quad q\geq \frac{N+\tau_+(\mu_2)}{-\tau_+(\mu_1)}; \qquad \\[6mm]
    (ii)  \qquad \ \ \, \frac{2-\tau_+(\mu_2)}{-\tau_+(\mu_1)}< q< \frac{N+\tau_+(\mu_2)}{-\tau_+(\mu_1)} \quad{\rm and}\quad \tau_+(\mu_1)(pq-1)+2p+2<0; \qquad\qquad \\[6mm]
     (iii)  \qquad  \ \ \frac{2-\tau_+(\mu_1)}{-\tau_+(\mu_2)} < p<\frac{N+\tau_+(\mu_1)}{-\tau_+(\mu_2)} \quad{\rm and}\displaystyle \quad \tau_+(\mu_2)(pq-1)+2q+2<0,
  \end{array}
$$
 then system (\ref{s 1.1}) admits no positive supersolutions.

\end{theorem}

\begin{remark}
The region of $(p,\,q)$ formed by Theorem \ref{teo 2} $(ii$-$iii)$
contains $A_1\cup A_2$, where

\vskip0.1cm
\noindent $A_1=\left\{(p,\,q):\,  \frac{2-\tau_+(\mu_1)}{-\tau_+(\mu_2)}\leq p< \frac{N+\tau_+(\mu_1)}{-\tau_+(\mu_2)} \quad {\rm and}\quad \frac{2-\tau_+(\mu_2)}{-\tau_+(\mu_1)}<  q< \frac{N+\tau_+(\mu_2)}{-\tau_+(\mu_1)}\right\}$

\vskip0.1cm
\noindent  and

\vskip0.1cm
\noindent
 $A_2=  \left\{(p,\,q):\,  \frac{2-\tau_+(\mu_1)}{-\tau_+(\mu_2)}< p< \frac{N+\tau_+(\mu_1)}{-\tau_+(\mu_2)} \quad {\rm and}\quad \frac{2-\tau_+(\mu_2)}{-\tau_+(\mu_1)}\leq  q< \frac{N+\tau_+(\mu_2)}{-\tau_+(\mu_1)}\right\}.
 $\end{remark}

 \begin{figure}[H]
    \centering
    \subfigure{
    \begin{minipage}{60mm}
  \includegraphics[scale=0.50]{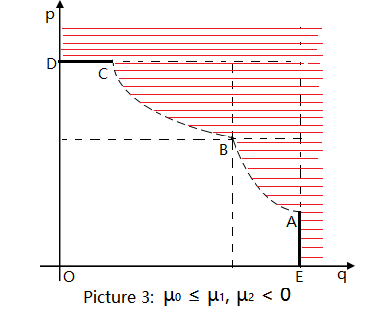}

  \end{minipage}
  }

  {\footnotesize  In Picture 3,  $E=(\frac{N+\tau_+(\mu_2)}{-\tau_+(\mu_1)},0)$,  $D=(0, \frac{N+\tau_+(\mu_1)}{-\tau_+(\mu_2)})$ and $B=(\frac{2-\tau_+(\mu_2)}{-\tau_+(\mu_1)},\, \frac{2-\tau_+(\mu_1)}{-\tau_+(\mu_2)})$.         The red regions of $(p,q)$ represent the non-existence showed in Theorem \ref{teo 2},  the dotted lines in the boundary of  the red regions are still open. }
 \end{figure}

Our method for the nonexistence of system (\ref{s 1.1}), different from the one in \cite{BDT,D},  is based on the nonexistence for nonhomogeneous Hardy problem
 \begin{equation}\label{eq 2.1}
 \arraycolsep=1pt\left\{
\begin{array}{lll}
 \displaystyle   \mathcal{L}_\mu u= f\quad
   {\rm in}\ \ {\Omega}\setminus \{0\},\\[1.5mm]
 \phantom{   L_\mu   }
 \displaystyle  u= 0\quad  {\rm   on}\ \, \partial{\Omega},
 \end{array}\right.
\end{equation}
when nonnegative function $f$ is locally H\"older continuous in $\bar\Omega\setminus\{0\}$
but $f\not\in L^1(\Omega,\, \Gamma_\mu dx)$, see reference \cite{CQZ1}.   By contradiction,  if there is a positive supersolution of system (\ref{s 1.1}), we  obtain  initial singularities at the origin, and then a contradiction could be obtained that one of the nonlinearities is not in related $L^1$ spaces; if not, we iterate these singularities to improve  new singularities until the nonlinearity is not in related $L^1$ spaces.

 \begin{theorem}\label{teo 3}
 $(i)$ Assume that  $\Omega=B_1(0)$, $\mu_0\leq \mu_1<0\leq \mu_2$ and $p,q>0$ satisfy that
\[
q< \frac{N+\tau_+(\mu_2)}{-\tau_+(\mu_1)} \quad {\rm and} \quad  \tau_+(\mu_1) (pq-1)+2p+2> 0,   \]
system (\ref{s 1.1}) has at least a positive supersolution.\smallskip

$(ii)$ Assume that $\mu_0\leq \mu_1,\,\mu_2<0$ and $p,q>1$ satisfy one of the following assumptions:
$$
 \arraycolsep=1pt
\begin{array}{lll}
   (a.1) \qquad  \,\frac{2-\tau_+(\mu_2)}{-\tau_+(\mu_1)} <q< \frac{N+\tau_+(\mu_2)}{-\tau_+(\mu_1)}\quad&{\rm and}\quad \tau_+(\mu_1)(pq-1)+2p+2>0;
   \qquad \\[6mm]
    (a.2) \qquad  \,q\leq \frac{2-\tau_+(\mu_2)}{-\tau_+(\mu_1)}\quad&{\rm and}\quad  p<\frac{N+\tau_+(\mu_1)}{-\tau_+(\mu_2)};
   \qquad \\[6mm]
    (b.1)  \qquad  \, \frac{2-\tau_+(\mu_1)}{-\tau_+(\mu_2)}<p<\frac{N+\tau_+(\mu_1)}{-\tau_+(\mu_2)}\quad &{\rm and}\ \quad \tau_+(\mu_2)(pq-1)+2q+2>0; \qquad\qquad  \\[6mm]
    (b.2)  \qquad  \, p\leq \frac{2-\tau_+(\mu_1)}{-\tau_+(\mu_2)}\quad &{\rm and}\ \ \quad q<\frac{N+\tau_+(\mu_2)}{-\tau_+(\mu_1)},  \qquad\qquad
    \end{array}
$$
system (\ref{s 1.1}) has at least a positive super solution.
\end{theorem}

It is worth noting that when $\mu_0=\mu_1<0\leq \mu_2$, we construct positive solution of system (\ref{s 1.1})
by constructing suitable super and sub solutions, see Picture 4 and Picture 5. Together with the nonexistence results,
we conclude that our conditions for $(p,q)$ is sharp.   However,  there are some critical cases are still open to determine the existence:

Picture 4: $\mu_0<\mu_1<0\leq \mu_2$, the critical curve $AQ:$
{\small $$\Big\{(p,q):\, \frac{2}{-\tau_+(\mu_1)}<q< \frac{N+\tau_+(\mu_2)}{-\tau_+(\mu_1)} \quad {\rm and} \quad  \tau_+(\mu_1) (pq-1)+2p+2=0\Big\}. $$}

Picture 5:   $\mu_0\leq \mu_1, \mu_2<0$,  point   $B=(\frac{2-\tau_+(\mu_2)}{-\tau_+(\mu_1)},\, \frac{2-\tau_+(\mu_1)}{-\tau_+(\mu_2)})$,   the critical curve $AB$:
{\small  $$\Big\{(p,q):\,  \frac{2-\tau_+(\mu_2)}{-\tau_+(\mu_1)} <q< \frac{N+\tau_+(\mu_2)}{-\tau_+(\mu_1)}\quad{\rm and}\quad \tau_+(\mu_1)(pq-1)+2p+2=0\Big\}$$}
 and $BC$
 {\small  $$\Big\{(p,q):\,  \frac{2-\tau_+(\mu_1)}{-\tau_+(\mu_2)} <p< \frac{N+\tau_+(\mu_1)}{-\tau_+(\mu_2)}\quad{\rm and}\quad \tau_+(\mu_2)(pq-1)+2q+2=0\Big\}.$$}

 \begin{figure}[H]
    \centering
   \subfigure{
  \begin{minipage}{65mm}
   \includegraphics[scale=0.50]{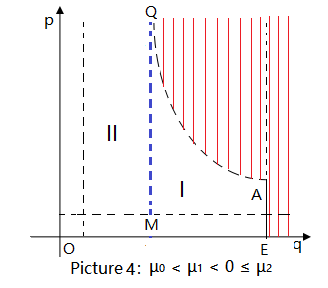}

  \end{minipage}
  }
    \subfigure{
    \begin{minipage}{60mm}
  \includegraphics[scale=0.47]{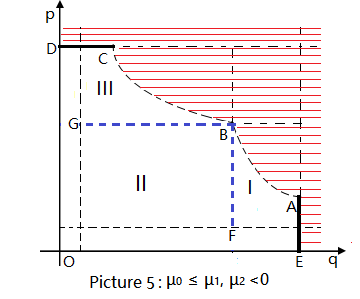}

  \end{minipage}
  }

  {\footnotesize    In regions I, II of Picture 4 and in regions I, II III of Picture 5, we construct the existence of positive solution of system (\ref{s 1.1})  in Theorem \ref{teo 3}.    }
 \end{figure} \smallskip

  The rest of this paper is organized as follows. In Section 2, we introduce the classification of single nonhomogeneous Hardy problem  and obtain the estimates of asymptotic behavior at the origin for specific nonhomogeneous term.   Section 3 is devoted to proving the nonexistence of  positive supersolutions by iterating to improve  original singularities to new ones  until they are not be in related spaces.
Section 4 is addressed to construct positive supersolutions in the subcritical case.

\setcounter{equation}{0}
\section{Preliminary results}
 Let us start this section with the Hardy problem (\ref{eq 2.1}). For $\mu_0 \le \mu <0$,
  $ \Phi_\mu$ is a regular solution of
(\ref{eq 1.1}) in the sense that $\mu |\cdot|^{-2} \Phi_{\mu}(\cdot) \in L^1_{loc}(\R^N)$ and
$$%\begin{equation}\label{eq 1.100}
\mathcal{L}_\mu \Phi_{\mu}= 0\quad {\rm in}\quad  \mathcal D'(\R^N),
$$%\end{equation}
that is,
$$  \int_{\R^N} \Phi_\mu \mathcal{L}_\mu \xi\, dx=0, \quad\forall\,\xi\in C^\infty_c(\R^N).$$
This normally distributional sense fails to express the singularity of $\Phi_\mu$. Furthermore, it even does not hold for $\Phi_\mu$ by the strong singularity for $\mu>0$ large.

  To overcome this difficulty, the authors  in \cite{CQZ1,CV1} introduced a new distributional framework: function $\Phi_\mu $ is  a  distributional solution of
$$%\begin{equation}\label{eq003}
\mathcal{L}_\mu   u =c_{\mu}\delta_0 \quad   {\rm in}\quad \mathcal D'(\R^N),
$$%\end{equation}
if $$%\begin{equation}\label{eq0030}
\int_{\R^N}\Phi_\mu \mathcal{L}_\mu^*\xi\, d\gamma_\mu =c_\mu\xi(0),\quad\forall\, \xi\in C_c^\infty(\R^N)
$$%\end{equation}
where $\delta_0$ is Dirac mass at the origin,
  $d\gamma_\mu(x) =\Gamma_\mu(x) dx$,
$$%\begin{equation}\label{L}
\mathcal{L}^*_\mu=-\Delta -2\frac{\tau_+(\mu) }{|x|^2}\,x\cdot\nabla,
$$%\end{equation}
and $c_{\mu}>0$ is the normalized constant.
With the help of this new distributional framework, the authors also obtain the classification of isolated singular
 solutions for non-homogeneous Hardy problem (\ref{eq 2.1}).
 \begin{theorem}\label{teo 2.1}\cite[Theorem 1.3]{CQZ1}
Let    $f$ be a function in $C^\gamma_{loc}(\overline{\Omega}\setminus \{0\})$ for some $\gamma\in(0,1)$.

 (I)  Assume that
\begin{equation}\label{f1}
\int_{\Omega} |f|\,   d\gamma_\mu <+\infty,
\end{equation}
  then    problem (\ref{eq 2.1}), subject to $\displaystyle \lim_{|x|\to0}u(x)\Phi_\mu^{-1}(x)=k$ with $k\in\R$,
has a  unique solution $u_k$, which  satisfies the distributional identity
 \begin{equation}\label{1.2f}
 \int_{\Omega}u_k  \mathcal{L}_\mu^*(\xi)\, d\gamma_\mu    = \int_{\Omega} f  \xi\, d\gamma_\mu  +c_\mu k\xi(0),\quad\forall\, \xi\in   C^{1.1}_0(\Omega).
 \end{equation}

 (II)  Assume that $f$ verifies (\ref{f1})  and $u$ is a nonnegative solution of (\ref{eq 2.1}), then  $u$ satisfies (\ref{1.2f}) for some $k\ge0$
and verifies that
$$\lim_{|x|\to0}u(x)\Phi_\mu^{-1}(x)=k.$$\smallskip

 (III)  Assume that $f\ge0$ and
\begin{equation}\label{f2}
 \lim_{r\to0^+} \int_{\Omega\setminus B_r(0)} f\, d \gamma_\mu   =+\infty,
\end{equation}
then problem (\ref{eq 2.1}) has no nonnegative solution.
\end{theorem}

In order to obtain the nonexistence of positive supersolution of system (\ref{s 1.1}),  we need the following basic tools:   comparison principle and basic estimate of  singularities at the origin.

\begin{lemma}\label{cr hp}
Assume that $\mu\geq \mu_0$, $f_1\geq f_2$, $k_1\geq k_2$,  functions $u_1$ and $u_2$ satisfy that $$%\begin{equation}\label{eq0 2.1}
 \arraycolsep=1pt\left\{
\begin{array}{lll}
 \displaystyle \mathcal{L}_\mu u_1 \geq f_1\qquad
   {\rm in}\quad  {\Omega}\setminus \{0\},\\[2mm]
  \phantom{   }
  \displaystyle \liminf_{|x|\to0}u_1(x)\Phi_\mu^{-1}(x)\geq k_1
 \end{array}\right.
$$%\end{equation}
and
$$%\begin{equation}\label{eq0 2.1}
 \arraycolsep=1pt\left\{
\begin{array}{lll}
 \displaystyle \mathcal{L}_\mu u_2 \leq f_2\qquad
   {\rm in}\quad  {\Omega}\setminus \{0\},\\[2mm]
 \phantom{   }
  \displaystyle \limsup_{|x|\to0}u_2(x)\Phi_\mu^{-1}(x)\leq k_2.
 \end{array}\right.
$$
  If $u_1\ge u_2$ on $\partial\Omega$, then
$$u_1\ge u_2\quad{\rm in}\quad \Omega\setminus \{0\}.$$
\end{lemma}
{\bf Proof.} Under the assumptions, $w:=u_1-u_2$ satisfies that
$$%\begin{equation}\label{eq0 2.1}
 \arraycolsep=1pt\left\{
\begin{array}{lll}
 \displaystyle \mathcal{L}_\mu w \geq 0\qquad
   {\rm in}\quad  {\Omega}\setminus \{0\},\\[1.5mm]
 \phantom{ \mathcal{L}_\mu   }
 w\geq 0\qquad
   {\rm on}\quad  {\partial\Omega}\\[1.5mm]
 \phantom{  }
  \displaystyle \liminf_{|x|\to0}w(x)\Phi_\mu^{-1}(x)\geq 0.
 \end{array}\right.
$$
 By applying \cite[Lemma 2.2]{CQZ1}, we obtain that $w\geq0$ and the proof ends.\hfill$\Box$

\begin{lemma}\label{cr 2.1}
Assume that  $\mu\geq\mu_0$, $r_0\in(0,1]$  is such that  $B_{2r_0}(0)\subset \Omega$ and nonnegative  function $u\in C^{1.1}(\Omega\setminus \{0\})$  is a classical solution of (\ref{eq 2.1})  with $f\gneqq  0$.

Then $u$ is positive in $\Omega\setminus\{0\}$ and for some $c_1>0$,
\begin{equation}\label{2.3-3}
u(x) \geq c_1 |x|^{\tau_+(\mu)}, \quad\forall\,x\in B_{r_0}(0)\setminus \{0\}.
\end{equation}
\end{lemma}
{\bf Proof.} Note that $u$ is positive in $\Omega\setminus\{0\}$ by strong Maximum Principle,   then
$u(x)\geq c_1$ on $\partial B_{r_0}(0)$ and
$$\liminf_{|x|\to0} u(x)\Phi_\mu(x)^{-1}\geq0.$$
By the compactness of $\partial B_{r_0}(0)$,   there exist $c_2>0$ such that
$$u \geq c_2\Gamma_\mu\quad{\rm on}\ \,\partial B_{r_0}(0),$$
thus by Lemma \ref{cr hp}, we have that $u\geq c_2\Gamma_\mu$ in $B_{r_0}(0)\setminus\{0\}$, that is (\ref{2.3-3}).  \hfill$\Box$

\begin{lemma}\label{lm 2.1-singular}
Assume that  $\mu\geq\mu_0$,   $f$ is a nonnegative continuous  function in $\Omega\setminus \{0\}$ and verifies that
\begin{equation}\label{2.1-f}
  f(x)\geq c_3 |x|^{\tau-2},\quad   \forall\, x\in B_{r_0}(0)\setminus\{0\},
\end{equation}
where   $\tau\in\R$ and $c_3\geq 1$.

 If $\tau\leq \tau_-(\mu)$, then
\begin{equation}\label{2.2-f}
 \arraycolsep=1pt\left\{
\begin{array}{lll}
 \displaystyle  \mathcal{L}_\mu u= f\qquad
   {\rm in}\quad \ {\Omega}\setminus \{0\},\\[1mm]
 \phantom{  L_\mu \, }
 \displaystyle  u= 0\qquad  {\rm   on}\quad\ \partial{\Omega},\\[1mm]
 \phantom{   }
  \displaystyle \lim_{|x|\to0}u(x)\Phi_\mu^{-1}(x)=0
 \end{array}\right.
\end{equation}
admits no positive solution

If $\tau> \tau_-(\mu)$ and $f\in L^1(\Omega,d\gamma_\mu)$, problem (\ref{2.2-f}) admits a unique solution  $u_f$, which is positive and has following asymptotic behavior:\\
for $\tau_-(\mu)<\tau<\tau_+(\mu)$,
\begin{equation}\label{2.3}
 u_f(x)\ge c_4|x|^{\tau},\ \    \forall\, x\in B_{r_0}(0)\setminus \{0\};
\end{equation}
 for $ \tau=\tau_+(\mu)$,
\begin{equation}\label{2.3-1}
 u_f(x)\ge c_4|x|^{\tau} (-\ln |x|),\ \     \forall\, x\in B_{r_0}(0)\setminus \{0\};
\end{equation}
and for $\tau>\tau_+(\mu)$,
\begin{equation}\label{2.3-2}
 u_f(x)\geq c_4|x|^{\tau_+(\mu)},\ \     \forall\, x\in  B_{r_0}(0)\setminus \{0\},
\end{equation}
where   $c_4\geq 1$.
\end{lemma}
{\bf Proof.} The nonexistence for $\tau\leq \tau_-(\mu)$ and existence for  $\tau>\tau_-(\mu)$ follow directly from Theorem \ref{teo 2.1}. For asymptotic behavior at the origin, we shall construct suitable  bounds for $u_f$. In the proof of this lemma, we assume that
$r_0=1$, otherwise we just do a scaling to re-normalized it.

For $\tau_-(\mu)<\tau< \tau_+(\mu)$,   we have that
  $$\mathcal{L}_\mu |x|^{\tau}=c_\tau |x|^{\tau-2},$$
where $c_\tau>0$. From (\ref{2.1-f}),  there exists $t_1>0$ such that
$$\mathcal{L}_\mu (t_1|x|^{\tau})\geq f(x),\ \, \forall\, x\in B_1(0)\setminus \{0\}.$$
Then by Lemma \ref{cr hp},    the lower bound is $t_1(|x|^{\tau}-|x|^{\tau_+(\mu)})$ in $B_1(0)\setminus\{0\}$.\smallskip

For  $\tau= \tau_+(\mu)$, direct computation shows that
$$\mathcal{L}_\mu (|x|^{\tau_+(\mu)}(-\ln |x|))=(2\tau_+(\mu)+N-2)|x|^{\tau_+(\mu)}. $$
The lower bound could be constructed by adjusting the coefficient $t$ of
$  t|x|^{\tau_+(\mu)}(-\ln |x|).$\smallskip

For  $\tau>\tau_+(\mu)$, we have that
$$\mathcal{L}_\mu (|x|^\tau)= (\mu-\tau(\tau+N-2))|x|^{\tau-2}, $$
where $\mu-\tau(\tau+N-2)<0$ for $\tau>\tau_+(\mu)$.
The lower bound could be constructed by adjusting the coefficients $t$ of
$ t(|x|^{\tau_+(\mu)}-|x|^{\tau}).$
The proof is thus completed. \hfill$\Box$

 \setcounter{equation}{0}
\section{Nonexistence in sharp supercritical regions}
\subsection{  Preliminary nonexistence property}
For $\mu\in[\mu_0,0)$, we note that $\Gamma_\mu$ defined in (\ref{1.2}) blows up at
the origin, which may cause  that  $u^q\not\in L^1(\Omega, d\gamma_{\mu_2})$ or
$v^p\not\in   L^1(\Omega, d\gamma_{\mu_1})$,  and then a contradiction appears by (\ref{f2})  in Theorem \ref{teo 2.1}.

\begin{proposition}\label{prop 2.1}
$(i)$  Let $\mu_0\leq \mu_1,\mu_2<0$,  and  $p,q>0$ satisfy
 \begin{equation}\label{qu 3}
 {\rm either}\qquad p\geq \frac{N+\tau_+(\mu_1)}{-\tau_+(\mu_2)}\quad{\rm or}\quad q\geq \frac{N+\tau_+(\mu_2)}{-\tau_+(\mu_1)},
 \end{equation}
 then system (\ref{s 1.1}) admits no positive supersolution.\smallskip\smallskip

 \noindent $(ii)$ Let $\mu_0\leq \mu_1<0\leq \mu_2$ and $p,q$ satisfy  one of the following conditions:
 \begin{equation}\label{qu 3}
 \arraycolsep=1pt
\begin{array}{lll}
   (a) \qquad  q\geq \frac{N+\tau_+(\mu_2)}{-\tau_+(\mu_1)};  \\[6mm]
   (b)  \qquad    \frac{2}{-\tau_+(\mu_1)}<q<\frac{N+\tau_+(\mu_2)}{-\tau_+(\mu_1)} \quad{\rm and}\quad  \tau_+(\mu_1)(pq+1)+2p+N\leq 0,\qquad\qquad\qquad
  \end{array}
 \end{equation}
 then system (\ref{s 1.1}) admits no positive supersolution.
\end{proposition}
{\bf Proof. }   By contradiction, we assume that $(u,v)$ is a positive solutions' pair of  system (\ref{s 1.1}), then by Lemma \ref{cr 2.1},   $u,v$ are positive and
\begin{equation}\label{2.4}
u(x)\geq c_5 |x|^{\tau_+(\mu_1)}\ \ {\rm and}\ \  v(x)\geq c_5 |x|^{\tau_+(\mu_2)}, \quad\forall\, x\in B_{r_0}(0)\setminus\{0\},
\end{equation}
thus,
$$%\begin{equation}\label{2.4-1}
u(x)^q\geq c_5 |x|^{\tau_+(\mu_1)q}\ \ {\rm and}\ \  v(x)^p\geq c_5 |x|^{\tau_+(\mu_2)p}, \quad \forall\, x\in B_{r_0}(0)\setminus\{0\},
$$%\end{equation}
where we recall that $r_0\in (0,1]$ satisfies that  $B_{2r_0}(0)\subset \Omega$.\smallskip

$(i)$ For $\mu_0\leq \mu_1,\mu_2<0$, we have that $\tau_+(\mu_1),\, \tau_+(\mu_2)<0$.

\smallskip

$\bullet$ If $\tau_+(\mu_2)p+\tau_+(\mu_1)\leq -N$, then $\mathcal{L}_{\mu_1}u\geq v^p\not\in L^1(\Omega,d\gamma_{\mu_1})$, thus,
 such function $u$ does not exist by Theorem \ref{teo 2.1} part (III).
 
\smallskip

$\bullet$ If $\tau_+(\mu_1)q+\tau_+(\mu_2)\leq -N$,
then we have that $\mathcal{L}_{\mu_2} v\geq u^q\not\in L^1(\Omega,d\gamma_{\mu_2}) $ and such function $v$ does not
 exist.
 
 \smallskip

$(ii)$ For $\mu_0\leq \mu_1<0\leq \mu_2$,  there holds that $\tau_+(\mu_1)<0< \tau_+(\mu_2)$.

\smallskip

$\bullet$ If $\tau_+(\mu_1)q\leq -(N+\tau_+(\mu_2))$, then $u^q\not\in L^1(\Omega, d\gamma_{\mu_2})$ and we derive  the nonexistence.

\smallskip

$\bullet$ If $-(N+\tau_+(\mu_2))< \tau_+(\mu_1)q< -2$, from the inequality
$$\mathcal{L}_{\mu_2} v(x)\geq u(x)^q \geq c_5^q|x|^{\tau_+(\mu_1)q},$$
 Lemma \ref{cr hp} and Lemma \ref{lm 2.1-singular} it follows that
$$v(x)\geq c_6|x|^{\tau_+(\mu_1)q+2}\quad{\rm for }\quad 0<|x|<r_0.$$
Then we have that $v(x)^p\geq c_6^p|x|^{(\tau_+(\mu_1)q+2)p}$ with
  $(\tau_+(\mu_1)q+2)p+\tau_+(\mu_1)\leq -N$.  Consequently, $\mathcal{L}_{\mu_1} u\not\in L^1(\Omega, d\gamma_{\mu_1})$ and we obtain the nonexistence from Theorem \ref{teo 2.1} part  (III).
\hfill$\Box$

\subsection{Nonexistence by iteration method }

\begin{proposition}\label{prop 3.1}
 Assume that $\mu_0\leq \mu_1<0\leq \mu_2$ and $p,q>0$ satisfy  that
 \begin{equation}\label{z 3.1}
    \frac{2}{-\tau_+(\mu_1)}  <q<   \frac{N+\tau_+(\mu_2)}{-\tau_+(\mu_1)}
      \quad{\rm
  and}\quad
        \tau_+(\mu_1) (pq-1)+2p+2< 0,
  \end{equation}
   then system (\ref{s 1.1}) admits no positive supersolution.
\end{proposition}
{\bf Proof.}
We argue by contradiction and let $(u,v)$ be a positive solutions' pair of system (\ref{s 1.1}).
Denote
\begin{equation}\label{d 3.1}
\tau_1^{(0)}=\tau_+(\mu_1)\quad{\rm and}\quad \tau_2^{(0)}=\tau_+(\mu_2).
\end{equation}
Then $u,v$ are positive and satisfy (\ref{2.4}), which could be written as
\begin{equation}\label{3.1}
u(x)\geq c_6 |x|^{\tau_1^{(0)}}\ \ {\rm and}\ \  v(x)\geq c_6 |x|^{\tau_2^{(0)}}, \ \ \, \forall\, x\in B_{r_0}(0)\setminus\{0\}.
\end{equation}

%\noindent {\it Subcritical case: $\tau_+(\mu_1) (pq-1)+2p+2<0$. }
Note that (\ref{z 3.1}) implies that  $\tau_-(\mu_2)<\tau_1^{(0)}q+2<0\leq \tau_+(\mu_2)$, then
$$\mathcal{L}_{\mu_2} v(x)\geq u(x)^q\geq c_6^q |x|^{\tau_1^{(0)}q}, \ \ \,\forall\, x\in B_{r_0}(0)\setminus\{0\},$$
which, by Lemma \ref{lm 2.1-singular}, implies that
$$v(x)\geq c_7|x|^{\tau_2^{(1)} }, \ \ \, \forall\, x\in B_{r_0}(0)\setminus\{0\},$$
where
$$\tau_2^{(1)}=\tau_1^{(0)}q+2\in\left(\tau_-(\mu_2),\,0\right).$$
Then by using the inequality
$$\mathcal{L}_{\mu_1} u(x)\geq v(x)^p\geq c_7^p |x|^{\tau_2^{(1)}p}, \ \ \, \forall\, x\in B_{r_0}(0)\setminus\{0\},$$
we have that
$$u(x)\geq c_8|x|^{\tau_1^{(1)} }, \ \ \, \forall\, x\in B_{r_0}(0)\setminus\{0\},$$
where   the second inequality of (\ref{z 3.1}) implies that
\begin{eqnarray*}
\tau_1^{(1)} & = &\tau_2^{(1)}p+2\\
                & = & (\tau_1^{(0)}q+2)p+2<\tau_1^{(0)}.
\end{eqnarray*}

If $\tau_1^{(1)}\leq \tau_-(\mu_1)$, then a contradiction could be obtained
from Theorem \ref{teo 2.1} part (III). If not, we use $\tau_1^{(1)}$ as the initial data, and
repeat above analysis.

In fact, we assume by induction that for $j=2,3,\cdots$
$$  v(x)\geq c_{j-1,1}|x|^{\tau_2^{(j-1)} }, \ \,  u(x)\geq c_{j-1,2}|x|^{\tau_1^{(j-1)} }, \ \ \forall\, x\in B_{r_0}(0)\setminus\{0\},$$
where
$$\tau_1^{(j-1)}\in(\tau_-(\mu_1),\, \tau_+(\mu_1)) \quad {\rm and}\quad \tau_2^{(j-1)}\in(\tau_-(\mu_2),\, \tau_+(\mu_2)) .$$
Repeating above analysis, we first obtain that
$v(x)\geq c_{j,1}|x|^{\tau_2^{(j)} }, \ \, \forall\, x\in B_{r_0}(0)\setminus\{0\}$
with
$$\tau_2^{(j)}=\tau_1^{(j-1)}q+2\in\left(\tau_-(\mu_2),\,\tau_2^{(0)}\right).$$
If $\tau_2^{(j)}\leq \tau_-(\mu_2)$, a contradiction is obtained and we are done.
If not, we continue to derive
$u(x)\geq c_{j,2} |x|^{\tau_1^{(j)} }, \ \, \forall\, x\in B_{r_0}(0)\setminus\{0\}$
with
$$\tau_1^{(j)}  = \tau_2^{(j)}p+2.$$
If $\tau_1^{(j)}\leq \tau_-(\mu_1)$,  we are done.

Note that
\begin{eqnarray*}
\tau_1^{(j)}   =  \tau_2^{(j)}p+2 =   (\tau_1^{(j-1)}q+2)p+2
\end{eqnarray*}
and
\begin{eqnarray*}
\tau_2^{(j)}   =  \tau_1^{(j)}q+2 =   (\tau_2^{(j-1)}p+2)q+2.
\end{eqnarray*}
{\bf Claim 1:} Assume that $pq>1$ and there exist $j_0\geq 1$ such that
\begin{equation}\label{start 1}
\tau_1^{(j_0)}< \tau_1^{(j_0-1)}\quad {\rm and}\quad \tau_2^{(j_0)}< \tau_2^{(j_0-1)},\end{equation}
then we have that
\begin{equation}\label{3.01}
\tau_1^{(j)}\to-\infty\quad {\rm and}\quad \tau_2^{(j)}\to-\infty\quad\ {\rm as}\ \, j\to+\infty.
\end{equation}
 In fact, let
$$s_j=\tau_1^{(j)}-\tau_1^{(j-1)}\quad{\rm and}\quad  t_j=\tau_2^{(j)}-\tau_2^{(j-1)},$$
then we have
$$s_j=pq s_{j-1}=\cdots=(pq)^{j_0-1}s_{j_0}\quad{\rm and}\quad t_j=pq t_{j_0-1}=\cdots=(pq)^{j_0-1}t_{j_0}.$$
From the fact that $pq>1$, we deduce  (\ref{3.01}).\smallskip

 From (\ref{z 3.1}),  we deduce that $pq>1$ and the sequence $\{\tau_1^{(j)}\}_j$ starts from $j=0$, $\{\tau_2^{(j)}\}_j$
does from $j=1$ and these two sequences are strictly decreasing. From {\bf Claim 1} there exists an integer $j_1\geq j_0$ such that
 $$\tau_1^{(j_1)}\leq \tau_-(\mu_1)\quad {\rm or}\quad \tau_2^{(j_1)}\leq \tau_-(\mu_2).$$
Thus, the nonexistence follows   from Theorem \ref{teo 2.1} part (III).\hfill$\Box$\smallskip

\begin{remark}
$(i)$ When $\mu_0<\mu_1<0\leq \mu_2$, the critical curve of $(p,q)$ is that
 $$q\in\left(\frac{2}{-\tau_+(\mu_1)},\,    \frac{N+\tau_+(\mu_2)}{-\tau_+(\mu_1)}\right)\quad{\rm and}\quad   \tau_+(\mu_1) (pq-1)+2p+2=0.$$ In this critical case,  we recall that
$$v(x)\geq c_9|x|^{\tau_2^{(1)} }, \ \, u(x)\geq c_9|x|^{\tau_1^{(1)} },\ \ \forall\, x\in B_{r_0}(0)\setminus\{0\},$$
where
$$\tau_2^{(1)}=\tau_1^{(0)}q+2\in\left(\tau_-(\mu_2),\,\tau_2^{(0)}\right)\quad{\rm and}\quad
\tau_1^{(1)} = \tau_2^{(1)}p+2 =\tau_1^{(0)}.$$
Observe that  the iteration procedure in Proposition \ref{prop 3.1} fails,  and  a new method has to be involved to obtain the nonexistence of super positive solution of system (\ref{s 1.1}). \smallskip

$(ii)$ According to the fact that $\tau_+(\mu_1)> -\frac{N-2}{2}$ for $\mu>\mu_0$,   (\ref{z 3.1})
is weaker than (\ref{qu 3}) case $(b)$. \\
For $\mu_1=\mu_0$,
the inequality
\begin{equation}\label{4.1}
\tau_+(\mu_1)(pq+1)+2p+N =\tau_+(\mu_1)(pq-1)+2p+2,
\end{equation}
We will take (\ref{qu 3}) as our assumption for the nonexistence .
\end{remark}

\begin{proposition}\label{prop 3.2}
 Assume that $\mu_0\leq \mu_1,\,\mu_2<0$ and $p,q>0$ satisfy  one of the following hypotheses:
$$%\begin{equation}\label{asum 3.1}
 \arraycolsep=1pt
\begin{array}{lll}
  (i) \qquad \ \ \frac{2-\tau_+(\mu_2)}{-\tau_+(\mu_1)}\leq q< \frac{N+\tau_+(\mu_2)}{-\tau_+(\mu_1)} \quad{\rm and}\quad \tau_+(\mu_1)(pq-1)+2p+2<0; \qquad\qquad \qquad\\[6mm]
    (ii)  \qquad  \ \frac{2-\tau_+(\mu_1)}{-\tau_+(\mu_2)} \leq p<\frac{N+\tau_+(\mu_1)}{-\tau_+(\mu_2)} \quad{\rm and}\displaystyle \quad \tau_+(\mu_2)(pq-1)+2q+2<0,
  \end{array}
$$%\end{equation}
 then system (\ref{s 1.1}) admits no positive supersolution.

\end{proposition}
{\bf Proof.}
We argue by contradiction and let $(u,v)$ be a positive supersolution  of  (\ref{s 1.1}).  Recall that $u,v$ are positive and satisfy (\ref{3.1}) with $\tau^{(0)}_1, \tau^{(0)}_2$  being defined in (\ref{d 3.1}).

 {\it Case $(i)$.}   From the inequality
 $$\mathcal{L}_{\mu_2} v(x)\geq u(x)^q\geq c_{10} |x|^{\tau_1^{(0)}q}, \ \ \, \forall\, x\in B_{r_0}(0)\setminus\{0\},$$
we derive  by Lemma \ref{lm 2.1-singular}   that
  $$v(x)\geq c_{11} |x|^{\tau^{(1)}_2}, \ \ \, \forall\, x\in B_{r_0}(0)\setminus\{0\},$$
  with
$$
 \tau_2^{(1)}=\min\{\tau_1^{(0)}q+2,\,\tau_2^{(0)}\} =\tau_1^{(0)}q+2 \in(\tau_-(\mu_2),\, \tau_2^{(0)}]
$$
from the assumption that $\frac{2-\tau_+(\mu_2)}{-\tau_+(\mu_1)}\leq q< \frac{N+\tau_+(\mu_2)}{-\tau_+(\mu_1)}.$

We use again the inequality
$$\mathcal{L}_{\mu_1} u(x)\geq v(x)^p\geq c_{11}^p |x|^{\tau_2^{(1)}p}, \quad \forall\, x\in B_{r_0}(0)\setminus\{0\}$$
to derive   that
  $$u(x)\geq c_{12} |x|^{\tau^{(1)}_1}, \quad  \forall\, x\in B_{r_0}(0)\setminus\{0\}$$
  with
  \begin{eqnarray*}
\tau_1^{(1)}  = \min\{\tau_2^{(1)}p+2,\,\tau_1^{(0)}\}
                & = & (\tau_1^{(0)}q+2)p+2<\tau_1^{(0)},
\end{eqnarray*}
by the assumption that
$$
\tau_+(\mu_1)(pq-1)+2p+2<0.
$$
If $\tau_1^{(1)}\leq \tau_-(\mu_1)$, then a contradiction could be obtained
from Theorem \ref{teo 2.1} part (III). If not, we use $\tau_1^{(1)}$ as the initial data, and
repeat above analysis.

We assume by induction that
 $$  v(x)\geq c_{j-1,1}|x|^{\tau_2^{(j-1)} }, \ \,  u(x)\geq c_{j-1,2}|x|^{\tau_1^{(j-1)} }, \ \ \forall\, x\in B_{r_0}(0)\setminus\{0\},$$
where
$$\tau_1^{(j-1)}\in(\tau_-(\mu_1),\, \tau_+(\mu_1)) \quad {\rm and}\quad \tau_2^{(j-1)}\in(\tau_-(\mu_2),\, \tau_+(\mu_2)) .$$
From above analysis, we derive that
$v(x)\geq c_{j,1}|x|^{\tau_2^{(j)} }, \ \, \forall\, x\in B_{r_0}(0)\setminus\{0\}$
with
$$\tau_2^{(j)}=\tau_1^{(j-1)}q+2\in\left(\tau_-(\mu_2),\,\tau_2^{(0)}\right).$$
If $\tau_2^{(j)}\leq \tau_-(\mu_2)$, a contradiction is obtained and we are done.
If not, we continue to deduce that
$u(x)\geq c_{j,2}|x|^{\tau_1^{(j)} }, \ \ \forall\, x\in B_{\frac12}(0)\setminus\{0\}$
with
$$\tau_1^{(j)}  = \tau_2^{(j)}p+2=\tau_1^{(j-1)}pq+2p+2.$$
 If $\tau_1^{(j)}\leq \tau_-(\mu_1)$,  we are done.
 Moreover, from above expression, we have that
 $$\tau_2^{(j+1)}  = \tau_1^{(j)}q+2=\tau_2^{(j)}pq+2q+2.$$
 If not, we start a new cycle above procedure and from {\bf Claim 1}, there exists $j_2\geq 2$ such that
 $$\tau_1^{(j_2)} \leq  \tau_-(\mu_1)\quad{\rm or}\quad \tau_1^{(j_2)} \leq  \tau_-(\mu_1)$$
 Therefore, the nonexistence of system (\ref{s 1.1}) follows.\smallskip

  {\it Case $(ii)$.} The nonexistence could be obtain by a similar proof.  \hfill$\Box$

\medskip

\noindent{\bf Proof of Theorem \ref{teo 1} and Theorem \ref{teo 2}.}   Theorem \ref{teo 1} $(i)$ follows Proposition  \ref{prop 2.1}   $(ii)$,  and  Theorem \ref{teo 1} $(ii)$ does
Proposition \ref{prop 2.1}  and observation (\ref{4.1}) for $\mu_1=\mu_0$,
and  Proposition \ref{prop 3.1} for $\mu_1\in(\mu_0,\, 0)$.\smallskip

Theorem \ref{teo 2} follows Proposition  \ref{prop 2.1}   $(i)$ and Proposition \ref{prop 3.2} directly.\hfill$\Box$

\setcounter{equation}{0}
\section{Existence in subcritical case }

We   construct super positive solutions of system (\ref{s 1.1})   by considering specific radially symmetric functions.

\medskip

\noindent{\bf Proof of Theorem \ref{teo 3}. } {\it Part $(i)$: We   construct super positive solutions of system (\ref{s 1.1}) when $\mu_0\leq \mu_1<0\leq \mu_2$  and $p,q>1$. } \smallskip

{\it Case 1: $q\in \left(\frac{2}{-\tau_+(\mu_1)},\,\frac{N+\tau_+(\mu_2)}{-\tau_+(\mu_1)}\right)$ and $\tau_+(\mu_1) (pq-1)+2p+2> 0$ (see the region I in picture 4). }

 Denote
  $$\tau_2=\tau_+(\mu_1)q+2\quad{\rm and} \quad \tau_1=(\tau_+(\mu_1)q+2)p+2. $$
 We note that $\tau_-(\mu_2)<\tau_2<\tau_+(\mu_2)$ and $\tau_1>\tau_+(\mu_1)$ by the assumptions that $\frac{2-\tau_+(\mu_2)}{-\tau_+(\mu_1)}< q<\frac{N+\tau_+(\mu_2)}{-\tau_+(\mu_1)}$ and $\tau_+(\mu_1) (pq-1)+2p+2> 0$, respectively.
 Let
$$u_1(x)=|x|^{\tau_+(\mu_1)}-|x|^{\tau_1}\quad{\rm and}\quad  v_1(x)=|x|^{\tau_2}\quad {\rm for}\ \ |x|>0.$$
Lemma \ref{lm 2.1-singular}  shows that for $c_{12}, c_{13}>0$,
 $$\mathcal{L}_{\mu_1} u_1(x)=c_{12} |x|^{(\tau_+(\mu_1)q+2)p}\geq c_{12} v_1(x)^p,\ \ \, \forall\, x\in B_1(0)\setminus\{0\}$$
 and
 $$\mathcal{L}_{\mu_2} v_1(x)= c_{13} |x|^{\tau_+(\mu_1)q }\geq  c_{13} u_1(x)^q,\ \ \, \forall\, x\in B_1(0)\setminus\{0\}.$$
Therefore, $t(u_1,v_1)$ is a positive solution of system (\ref{s 1.1}) in $B_1(0)\setminus\{0\}$ for $t\in(0,1)$ small enough.\smallskip

{\it Case 2: $q< \frac{2}{-\tau_+(\mu_1)} $ (see the region II in picture 4).}
Denote
  $$\tau_3= \tau_+(\mu_2)p+2 \quad{\rm and} \quad \tau_4=\tau_+(\mu_1)q+2, $$
and we have that   $\tau_3>0\geq \tau_+(\mu_1)$ and $\tau_4>\tau_+(\mu_2)\geq 0$.

 Let
$$u_2(x)=|x|^{\tau_+(\mu_1)}-|x|^{\tau_3}\quad{\rm and}\quad  v_2(x)=|x|^{\tau_+(\mu_2)}-|x|^{\tau_4}\quad {\rm for}\ \ |x|>0.$$
There exist  $c_{14},c_{15}>0$ such that
 $$\mathcal{L}_{\mu_1} u_2(x)=c_{14} |x|^{\tau_+(\mu_2) p}\geq c_{14} v_2(x)^p,\ \ \, \forall\, x\in B_1(0)\setminus\{0\}$$
 and
 $$\mathcal{L}_{\mu_2} v_2(x)= c_{16} |x|^{\tau_+(\mu_1)q }\geq  c_{15} u_2(x)^q,\ \ \, \forall\, x\in B_1(0)\setminus\{0\}.$$
Therefore, $t(u_2,v_2)$ is a positive solution of system (\ref{s 1.1}) in $B_1(0)\setminus\{0\}$ for $t\in(0,1)$ small enough.\smallskip

{\it Case 3: $q= \frac{2}{-\tau_+(\mu_1)} $ (line $MQ$ in picture 4).}
Denote     $\tau_5= \tau_+(\mu_2)p+1 $ and we have that   $\tau_5>0\geq \tau_+(\mu_1)$.
 Let
$$u_3(x)=|x|^{\tau_+(\mu_1)}-|x|^{\tau_5}\quad{\rm and}\quad  v_3(x)=|x|^{\tau_+(\mu_2)}(-\log|x|),\ \ \,\forall\, x\in B_1(0)\setminus\{0\}.$$
There exist  $r_1\in(0,r_0)$ and $c_{16},c_{17}>0$ such that
 $$\mathcal{L}_{\mu_1} u_3(x)=c_{16} |x|^{\tau_+(\mu_2) p-1}\geq c_{16}  v_3(x)^p,\ \ \, \forall\, x\in B_{r_1}(0)\setminus\{0\}$$
 and
 $$\mathcal{L}_{\mu_2} v_3(x)= c_{17} |x|^{\tau_+(\mu_2) -2}\geq  c_{17} u_3(x)^q,\ \ \, \forall\, x\in B_{r_1}(0)\setminus\{0\}.$$
Therefore, $t(u_3,\, v_3)$ is a positive solution of system (\ref{s 1.1}) in $B_{r_1}(0)\setminus\{0\}$ for $t\in(0,1)$ small enough. By scaling, a solution could be found in $B_1(0)\setminus\{0\}$.  \smallskip

\noindent{\it Part $(ii)$. We   construct super positive solutions of system (\ref{s 1.1}) when $\mu_0\leq \mu_1,\, \mu_2<0$ and $p,q>1$.}

\smallskip

{\it Case 4: $q\in \left(\frac{2-\tau_+(\mu_2)}{-\tau_+(\mu_1)},\,\frac{N+\tau_+(\mu_2)}{-\tau_+(\mu_1)}\right)$ and $\tau_+(\mu_1)(pq-1)+2p+2>0$ (see the region I in picture 5).}

The construction of a positive supersolution   follows
by functions' pair $t(u_1,\,v_1)$,  where $t>0$,
$$u_1(x)=|x|^{\tau_+(\mu_1)}-|x|^{\tau_1}\quad{\rm and}\quad  v_1(x)=|x|^{\tau_2}\quad {\rm for}\ \ |x|>0$$
with
  $$\tau_2=\tau_+(\mu_1)q+2\quad{\rm and} \quad \tau_1=(\tau_+(\mu_1)q+2)p+2. $$
  Note that
  $$\tau_-(\mu_2)<\tau_2<\tau_+(\mu_2)\quad{\rm and}\quad \tau_1>\tau_+(\mu_1).$$
 The proof is the same as  {\it Proof of Theorem \ref{teo 3}  $(i)$ case 1} and   we omit the details.
  \smallskip

{\it Case 5: $ q< \frac{2-\tau_+(\mu_2)}{-\tau_+(\mu_1)} $ and $ p<\frac{2-\tau_+(\mu_1)}{-\tau_+(\mu_2)}$  (see the region II in picture 5).}
We recall that
$$u_2(x)=|x|^{\tau_+(\mu_1)}-|x|^{\tau_3}\quad{\rm and}\quad  v_2(x)=|x|^{\tau_+(\mu_2)}-|x|^{\tau_4}\quad {\rm for}\ \ |x|>0$$
with   $$\tau_3= \tau_+(\mu_2)p+2 \quad{\rm and} \quad \tau_4=\tau_+(\mu_1)q+2.$$
Here $\tau_+(\mu_2)<\tau_4<0$ and $\tau_3>\tau_+(\mu_1)$. The rest of construction in
this case is routine and we omit the details.  \smallskip

{\it Case 6: $q=\frac{2-\tau_+(\mu_2)}{-\tau_+(\mu_1)} $ and $p<\frac{2-\tau_+(\mu_1)}{-\tau_+(\mu_2)}$ (see the line segment $BF$ in picture 5). }

 Let $\epsilon_0>0$ such that
$ \tau_+(\mu_2)p+2+\epsilon_0>\tau_+(\mu_1)$. Now we denote
$$u_4(x)=|x|^{\tau_+(\mu_1) }-|x|^{\tau_6} \quad{\rm and}\quad  v_4(x)=|x|^{\tau_+(\mu_2)}(-\ln|x|) \quad {\rm for}\ \ x\in B_1(0)\setminus\{0\},$$
where  $\tau_6=\tau_+(\mu_2)p+2+\epsilon_0>\tau_+(\mu_2)$ in this case. The rest  is routine and we omit the details.

 \smallskip

{\it  Case 7:  $0<q< \frac{2-\tau_+(\mu_2)}{-\tau_+(\mu_1)},\  p=\frac{2-\tau_+(\mu_1)}{-\tau_+(\mu_2)}$ (segment $BG$ in picture 5). }
A positive supersolution could be built by following functions:
$$u_6(x)=|x|^{\tau_+(\mu_1)}(-\ln|x|) \quad{\rm and}\quad  v_6(x)=|x|^{\tau_+(\mu_2)}-|x|^{\tau_8}\quad {\rm for}\ \ x\in B_1(0)\setminus\{0\},$$
where $\tau_8= \tau_+(\mu_1) q +2>\tau_+(\mu_2)$. \smallskip

{\it Case 8:  $p\in \left(\frac{2-\tau_+(\mu_1)}{-\tau_+(\mu_2)},\,\frac{N+\tau_+(\mu_1)}{-\tau_+(\mu_2)}\right)$ and $\tau_+(\mu_2)(pq-1)+2q+2>0$ (see the region III in picture 5).}

The construction of a positive supersolution   follows
by functions' pair $t(u_7,\,v_7)$,  where $t>0$,
$$v_7(x)=|x|^{\tau_+(\mu_2)}-|x|^{\tau_9}\quad{\rm and}\quad  u_7(x)=|x|^{\tau_{10}}\quad {\rm for}\ \ x\in B_1(0)\setminus\{0\}$$
with
  $$\tau_{10}=\tau_+(\mu_2)p+2\quad{\rm and} \quad \tau_9=(\tau_+(\mu_2)p+2)q+2. $$
  Note that
  $$\tau_-(\mu_1)<\tau_{10}<\tau_+(\mu_1)\quad{\rm and}\quad \tau_9>\tau_+(\mu_2).$$
 The rest is routine and we omit the details.  The proof is therefore complete.\hfill$\Box$

   \bigskip \bigskip

\noindent{\small {\bf Acknowledgements:} H. Chen is supported by NNSF of China, No: 12071189 and 12001252,
 by the Jiangxi Provincial Natural Science Foundation, No: 20202BAB201005,  20202ACBL201001.
B. Zhang is supported by the
National Natural Science Foundation of China (No. 11871199), the Heilongjiang Province Postdoctoral
Startup Foundation, P.R. China (LBH-Q18109), and the Cultivation Project of Young and Innovative Talents
in Universities of Shandong Province.

%\bibliographystyle{IEEEtran}%

%\bibliography{Habib}

\end{document}